\input amstex
\input amams.sty
\document
\annalsline{151}{2000}
\received{November 30, 1998}
\startingpage{269}
\font\tenrm=cmr10
\def\ritem#1{\item{{\rm #1}}}


\def\MD{\hbox{MD}}

\define\ti#1{\tilde #1}

\define\R{\Bbb R^N}
\define\N{\Bbb N}
\define\M{\Bbb M}

\define\na#1{|\nabla#1|^p}
\define\np#1{|\nabla#1|^{p-2}\nabla #1}
\define\ch#1{\chi_{\{#1 > 0\}}}

\define\J#1#2#3{\int_{\R} \left({\na#1\over p} -#2#1 + \left({p-1\over p}\right)#3^p\ch#1\right)dx}
\define\I#1#2{\int_{\R} \left({\na#1\over p}  + \left({p-1\over p}\right)#2^p\ch#1\right)dx}

\define\de#1#2{\frac{\partial#1 }{\partial #2 }}

\def\V{\Vert}

\catcode`\@=11
\font\twelvemsb=msbm10 scaled 1100

\font\ninemsb=msbm10 scaled 800
\newfam\msbfam
\textfont\msbfam=\twelvemsb  \scriptfont\msbfam=\ninemsb
  \scriptscriptfont\msbfam=\ninemsb
\def\msb@{\hexnumber@\msbfam}
\def\Bbb{\relax\ifmmode\let\next\Bbb@\else
 \def\next{\errmessage{Use \string\Bbb\space only in math
mode}}\fi\next}
\def\Bbb@#1{{\Bbb@@{#1}}}
\def\Bbb@@#1{\fam\msbfam#1}
\catcode`\@=12

 \catcode`\@=11
\font\twelveeuf=eufm10 scaled 1100
\font\teneuf=eufm10
\font\nineeuf=eufm7 scaled 1100
\newfam\euffam
\textfont\euffam=\twelveeuf  \scriptfont\euffam=\teneuf
  \scriptscriptfont\euffam=\nineeuf
\def\euf@{\hexnumber@\euffam}
\def\frak{\relax\ifmmode\let\next\frak@\else
 \def\next{\errmessage{Use \string\frak\space only in math
mode}}\fi\next}
\def\frak@#1{{\frak@@{#1}}}
\def\frak@@#1{\fam\euffam#1}
\catcode`\@=12
\title{Regularity of a free boundary with\\ application to the Pompeiu problem}
\shorttitle{Regularity of a free boundary}
\acknowledgements{L. Caffarelli was  partially supported  by the NSF.  H.\ Shahgholian was partially supported  by the
Swedish
 Natural Sciences Research Council, Wenner-Gren Center, and
 the Carl-Trygger Foundation.
He also thanks the University of Texas at Austin for a visiting
appointment.\hfill\break
\phantom{11}\hskip.25in 1991 {\it Mathematics Subject Classification}.  35R35, 31B20.
\hfill\break
\phantom{11}\hskip.25in {\it Key words and phrases}. Free boundary problems, regularity, Pompeiu problem.
}
 \twoauthors{Luis A. Caffarelli, Lavi Karp,}{Henrik Shahgholian }
   \institutions{University of Texas at Austin, 
Austin, TX\\
 {\eightpoint {\it E-mail address\/}: caffarel$\@$math.utexas.edu}\\
\vskip6pt
 Department of Applied Mathematics, International College of
Technology,\\  Karmiel,  Israel\\
{\eightpoint {\it E-mail address\/}:  karp$\@$techunix.technion.ac.il}
 \\
\vskip6pt
Royal Institute of Technology,  Stockholm, Sweden\\
{\eightpoint {\it E-mail address\/}: henriks$\@$math.kth.se}
 }

\bigbreak \centerline{\bf Abstract}
\medbreak In the unit ball $B(0,1)$, let $u$ and  $\Omega$ (a domain
  in $\R$) solve the following overdetermined problem:
$$\Delta u =\chi_\Omega\quad \hbox{ in } B(0,1), \qquad 0 \in \partial
\Omega, \qquad u=|\nabla u |=0 \quad \hbox{ in } B(0,1)\setminus \Omega,$$
where $\chi_\Omega$ denotes the characteristic function, and the equation is
satisfied in the sense of distributions.

If the complement of $\Omega$  does not develop cusp singularities
 at  the origin then we prove
$\partial \Omega$ is analytic in some  small neighborhood of the
origin.
The result can be modified to yield for more general divergence form 
operators. As an application of this, then, we obtain the regularity
of  the boundary  of a  domain without the   Pompeiu property,
provided  its complement has no cusp singularities.

\section{Introduction}

In this paper we study the regularity properties of solutions to
 a certain type of free boundary problems, resembling the obstacle problem
but with no sign assumption, i.e., with  no obstacle.
Mathematically the problem is formulated as follows.
Let $\Omega \subset \R$  and suppose there is a function $u$, solving
 the following overdetermined problem
$$\Delta u = \chi_\Omega \ \hbox{ in } B, \qquad \qquad u=|\nabla u|= 0 
\ \hbox{ in } B\setminus \Omega,\tag {1.1}$$
where $B$ is the unit ball. 

The question  we ask is whether $\partial \Omega$
is smooth. Indeed, if $\partial \Omega$ is an analytic surface then 
 by the  Cauchy-Kowalewski theorem, we can always solve
 the above overdetermined problem locally. We thus ask the
reverse of the Cauchy-Kowalewski theorem.

The problem also has  a
potential theoretic interpretation which is as follows.
Denote by $U$ the Newtonian potential of $\Omega$ (bounded set)
with constant density (i.e., the convolution of the fundamental solution
with $\chi_\Omega$) and with $0\in  \partial \Omega$. Suppose there
exists a harmonic function $w$ in $B(0,r)$ ($r$ small) such that $w=U$ in
$\Omega^c$ (the complement of $\Omega$); observe that   $U$ is harmonic in
$\Omega^c$.
 This property is called harmonic continuation.

 Next,  up to a normalization constant, $w-U$ satisfies equation (1.1).
Once again the question is whether possession of such a property,
harmonic
continuation for $\Omega$,
will result in the regularity of $\partial \Omega$ near the origin.
For the interested reader we refer to [Sa1--3] for similar types of problems.

It is noteworthy that this kind of problem has been in focus
of attention in mathematical physics, especially in geophysics [St], [Ma],
and  in inverse potential theory [I1].  Also, newly developed
problems in operator theory  reduce the analysis of the spectrum of
certain hyponormal operators to the study of the  
 solutions in problem (1.1), in the complex plane [MY], [P].
We also refer to an excellent book by H.\ S. Shapiro [Shap2], for
some basic aspects.

In order to state our main results, let us define a local solution.

\numbereddemo{{D}efinition} We say  a function $u$  belongs to the class
$P_r(z,M)$ if $u$ satisfies (in the sense of distributions):
\medbreak
\item{(1)} $\Delta u = \chi_\Omega \qquad \text{in } B_r(z)$,
\smallbreak\item{(2)} $u=|\nabla u|=0$ 
in $B_r(z)\setminus \Omega$,
\smallbreak\item{(3)} $\Vert u \Vert_{\infty,B_r(z)} \leq M$,
\smallbreak\item{(4)} $z \in \partial \Omega$.
\enddemo

We also denote by $P_\infty(0,M)$  ``global solutions''
 with quadratic growth, i.e., solutions 
in the entire space $\R$ with quadratic growth $|u(x)| \leq M (|x|^2 +
 1)$.

\demo{{R}emark} If $u \in P_r(z,M)$, then
\medbreak
\item{(1)}  $u(z+x) \in P_r(0,M)$,
\smallbreak\item{(2)} $u(z+rx)/r^2 \in P_1(0,M/r^2)$,
\smallbreak\item{(3)} Also if $\V D_{ij} u \V \leq M$, then  $u(z+rx)/r^2 \in P_1(0,M)$.
\enddemo
 
Obviously (1) implies that the class is point independent, so we can always
consider the class $P_r(0,M)$.
Our class $P$ differs from that of [Ca2] in two ways. We do not 
restrict the function $u$ to be nonnegative and we replace the uniform $C^{1,1}$-norm with a uniform $C^0$-norm.

This new feature introduces  new difficulties, and as a first task we have to
cope with the  optimal regularity of the function $u$ itself before attacking
the regularity problem for the boundary $\partial \Omega$.
The main tools in our study of problem (1.1) will be  a  monotonicity
lemma due to [ACF]; see Lemma 2.1 below.

\demo{Definition {\rm 1.2 (Minimal Diameter)}}
The minimum diameter of a  bounded set $D$, denoted $\MD(D)$, is the infimum of
distances between pairs
of parallel planes such that $D$ is contained in the strip determined by
the planes. We also define the density function
$$\delta_r (u) ={\MD (\{u=|\nabla u|=0\}\cap B(0,r)) \over r}.
$$
\enddemo
Now we can state our main results. In Section 7, we apply Theorems
I--III
to obtain the regularity of a domain without the Pompeiu property (see
\S 7
for details).

\nonumproclaim{Theorem I} There is a constant $C_1=C_1(N)$ such that if $u \in P_1(z,M)$
then 
$$\sup_{B(z,1/2)}\V D_{ij}u \V \leq C_1M.$$
\endproclaim

\nonumproclaim{Theorem II} Let $u \in P_\infty(0,M)$ and suppose $\{u =0\}$ has
nonempty interior{\rm ,} or $\delta_r(u) >0$ for some $r>0${\rm .} Then $ u\geq 0$
in $\R$ and   $D_{ee} u \geq 0$ in 
$\Omega${\rm ,} for any direction $e${\rm ,} i{\rm .}e{\rm .}\ $\Omega^c$ is convex{\rm .} Moreover
if
$\varlimsup_{R\to \infty} \delta_R(u) >0$ then $u$ is a half\/{\rm -}\/space
solution{\rm ,}
i{\rm .}e{\rm .,} $u=(\max(x_1,0))^2/2$ in some coordinate system{\rm .}
\endproclaim

\nonumproclaim{Theorem III} There exists a modulus of continuity $\sigma$
$(\sigma (0^+)=0)$ such that if $u \in P_1(0,M)$ and
$\delta_{r_0}(u) >\sigma(r_0)$ for some $r_0<1${\rm , }
then $\partial \Omega$ is the graph
of a $C^1$ function in $B(0,c_0r_0^2)${\rm .} Here $c_0$ is a universal
constant{\rm ,} depending only on $M$ and dimension{\rm .}
\endproclaim  

Subject to the condition in  Theorem III, the analyticity of the free
boundary now follows by classical results  [KN], [I2]. We thus have
the following corollary.

\nonumproclaim{Corollary to Theorem III} Under the thickness condition in
Theorem {\rm III,}  the free boundary in {\rm (1.1)} is
analytic{\rm ,} in some neighborhood of the origin{\rm .}
\endproclaim

Theorems I--III  are known  in 2-space dimensions [Sa2--3]. Indeed,
M. Sakai
[Sa2--3] gives a complete description of the boundaries 
of all such domains in ${\Bbb R}^2$.
Recently a different approach to this problem was made by
B. Gustafsson and M. Putinar [GP], where they proved
that $\partial \Omega$  in (1.1) is contained in an  analytic arc
in ${\Bbb R}^2$. This question is still open in higher dimensions.
For the case $u \geq 0$ (the original obstacle problem)
the first author has recently proved  the following:

\nonumproclaim{Structure of the singular set {\rm ([Ca4])}}
Let $u\in P_1(0,M)$ and $u \geq 0${\rm .}  Let 
$y$ also be a singular point of the free boundary in {\rm (1.1);} i{\rm .}e{\rm .,} the 
free boundary does not satisfy the condition in Theorem {\rm III,} near $y${\rm .}
Then there exists a unique nonnegative quadratic polynomial
{\rm (}\/and a unique matrix $A_{y}${\rm )}
$$Q_{y}=\frac12 x^T A_{y} x$$
 with {\rm trace}$A_{y}=1$ and  such that
\medbreak
\ritem{(1)} $|( u - Q_{y})(x)| \leq |x-y|^2\sigma(|x-y|),$
for some universal modulus of continuity $\sigma${\rm ,}
depending on $M$ only{\rm .}
\ritem{(2)}  $A_{y}$ is continuous on $y${\rm .}
\ritem{(3)} If ${\rm dim(kernel}(A_{y})) =k$ then there exists
 a $k$\/{\rm -}\/dimensional $C^1$ manifold ${\Cal T}_{y,u}${\rm ,} such that
 $$S_u \cap B(y,r) \subset {\Cal T}_{y,u}, $$
for some small $r${\rm .} Here $S_u$ indicates the singular points of the
free
boundary{\rm ,} i{\rm .}e{\rm .,} points which do not fall under the hypothesis of
Theorem {\rm III.}

\endproclaim
 
This fact with no positivity assumption,
is studied in a forthcoming paper
by  the first and the third authors.
In this paper, however,  we will prove the analyticity of the free
boundary only in the case
of ``thick'' complement described in Theorem III.

\demo{Plan of the paper} Section 2 is devoted to some
technical
tools, which are known, but probably not well known in the context
used in this paper. In Section 3 we carry out the proof for Theorem I,
using the monotonicity lemma (Lemma 2.1) and the blow-up technique.

In Section 4 we introduce further lemmas, which will  somehow exhaust
properties of the monotonicity lemma. These are used to prove Theorems
II--III
in Sections 5--6, respectively.
In Section 7 we generalize the monotonicity lemma to yield for
divergence type
operators 
$$\sum D_i ( a_{ij}D_j u) + a(x) u.$$
As a result we obtain the regularity of a domain without the Pompeiu
property, which  we explain to some extent in Section 7.
\enddemo

\section{The monotonicity formula}

In this section we will gather all basic tools  used to prove  
Theorems I--III.
A  fundamental tool, however, is the following monotonicity
 lemma.

\proclaimtitle{[ACF]}
\proclaim{Lemma} 
Let $h_1${\rm ,} $h_2$ be two nonnegative continuous sub\/{\rm -}\/solutions of $\Delta u =0$ in $B(x^0,R)$ ($R>0$){\rm .}
Assume further that $h_1h_2=0$ and that $h_1(x^0)=h_2(x^0)=0${\rm .} Then the following function
is monotone in $r$ {\rm (}$0<r<R${\rm )}
$$\varphi (r)={ 1 \over r^{4}}\left( \int_{B(x^0,r)} {| \nabla h_1|^2\over |x-x^0|^{N-2}}\right)\left( \int_{B(x^0,r)} {| \nabla 
h_2|^2 \over |x-x^0|^{N-2}}\right). \tag 2.1
$$
\endproclaim

For a fixed direction $e$, set 
$$(D_e u)^+ = \max(D_e u , 0) , \qquad \qquad (D_e u)^- = - \min(D_e u ,
0).$$
Then we have the following lemma.

\proclaim{Lemma} Let $u \in P_1(0,M)$ and consider the monotonicity formula
{\rm (2.1)} for the nonnegative subharmonic  functions $(D_e u)^\pm$  where $e$ is any fixed direction{\rm ,}
and denote this by $\varphi(r,D_e u)${\rm .}
Then 
$$\varphi'(r,D_e u) \geq {2 \over r} \varphi(r,D_e u)
 (\gamma(\Gamma_+) + \gamma(\Gamma_-) -2 ),
\tag 2.2$$
where 
$$\gamma(\Gamma_\pm) ( \gamma(\Gamma_\pm) + N -2) =
 \lambda(\Gamma_\pm),\tag {2.3}$$
with
$$\lambda(\Gamma_\pm)= 
\inf \left( { \int_{\Gamma_\pm}  |\nabla' w|^2 d\sigma \over
    \int_{\Gamma_\pm}  | w|^2 d\sigma  } \right).\tag 2.4$$
Here{\rm ,} $\nabla'$ is the gradient on $\Bbb S^{N-1}${\rm ,}
$\Gamma_{\pm}=\Gamma_{\pm}(r)$ is the projection of 
$ \partial B(0,r) \cap \{ (D_e u)^\pm >0\}$
onto the unit sphere{\rm ,} and the infimum has been taken over all
nonzero functions with compact support in  $\Gamma_\pm${\rm . }
\endproclaim

For a proof of this lemma see (the proof of) Lemma 5.1 in [ACF].
The ``set function'' $\gamma (E)$ as a function of 
$E \subset {\Bbb  S}^{N-1}$ is
  called   the characteristic constant of $E$; see, for example, the paper by Friedland and Hayman [FH].
A result of Sperner [Spn] states that, among all sets with given
  ($N-1$)-dimensional surface area $s \omega_N $ on the unit sphere
${\Bbb S}^{N-1}$ in $\R$, a spherical cap, i.e.\ a set of the form
$-1\leq c< x_1 \leq 1$, has the smallest characteristic constant
$\gamma(s,N)$, where $\omega_N$ is the area of the unit sphere in $\R$, 
$0<s<1$, and $c$ and $s$ are coupled by the relations
$$s=\frac{\omega_{N-1}}{\omega_N}\int_0^{\theta_0} (\sin t )^{N-2} dt,
\qquad c=\cos \theta_0, \qquad 0<\theta_0 < \pi.$$

It thus suffices to consider the function $\gamma$ as a function of
$s={\rm area}(E)/\omega_N$, with $s$ as above. 
From this [FH,  Thm.\  2]  deduce that for fixed $s$, $\gamma$ is a
monotone decreasing function in $N$ (the space dimension),
 and that the limit exists
as $N$ tends to infinity
$$\gamma (s,N) \geq \gamma(s,\infty).$$
On the other side as $N \to \infty$
we will  have, for some $h \geq -\infty$,
$$s =\frac{1}{\sqrt{2\pi} }\int_h^\infty \exp{(-t^2/2)} dt; \tag 2.5$$
 see e.g. [FH, p.~149]. 

Hence for any set $E \in {\Bbb S}^{N-1}$ and $s$ as above
$$\gamma(E,N)\geq\gamma(E^\star,N)= \gamma(s,N) \geq 
\lim_{N \to \infty} \gamma(s,N)=\gamma(s,\infty),\tag 2.6$$
where $E^\star$ is the above described symmetrization of $E$.

Now from (2.3) (cf.  [ACF, p.\  441]) one has
$$\gamma  =\frac{\lambda}{(N-2)} + O(N^{-3}).\tag 2.7$$
In (2.7) if we let $N $ tend to infinity we will  obtain,
by  some tedious   calculations (cf.  [FH, p.~149]), 
$$\frac{\lambda}{(N-2)} + O(N^{-3}) \to \Lambda,\tag 2.8$$
where $-\Lambda$ is the first Dirichlet eigenvalue of the
one-dimensional
Ornstein-Uhlenbeck operator $\Delta - x\cdot \nabla$ on the set
$(h,\infty)$, with $h$ as in (2.5).
Now (2.2)--(2.8) imply that

$$\varphi'(r,D_e u) \geq {2 \over r} \varphi(r,D_e u)
 (\Lambda(h_+) + \Lambda(h_-) -2 ),
\tag 2.9$$
where $h_\pm$ are the corresponding constants, in (2.5),
 for $\Gamma_\pm$.
Obviously $h_+ + h_- \geq 0$. Also, by  results of
Beckner-Kenig-Pipher [BKP] (cf. also [CK,  \S  2.4])
$\Lambda$ is convex  and $\Lambda(0) = 1$; hence 
$$\Lambda(h_+) + \Lambda(h_-) \geq 2\Lambda(a)= 2.\tag 2.10$$
In [BKP] it is actually proved that the convex function
$\Lambda$ satisfies
$$\Lambda''(0)=4(1-\ln 2)/\pi >0.$$

 For convenience  we now set
$$\gamma(r) = \gamma(\Gamma_+(r)) + \gamma(\Gamma_-(r)) -2.$$
Then inserting (2.10) in (2.9)  we  obtain the following lemma.

\proclaim{ Lemma} There holds $\gamma(r) \geq 0$ for all $r${\rm .}
 Moreover the strict inequality holds
unless $\Gamma^\star_\pm(r)$ are both half\/{\rm -}\/spheres{\rm .} In particular if any of
 the $\Gamma^\star_\pm (r)$ digresses from being a half\/{\rm -}\/spherical cap by an
 area\/{\rm -}\/size of $\varepsilon${\rm ,} say{\rm ,}
then 
$$\gamma(r) \geq C\varepsilon^2 .\tag 2.11$$
\endproclaim

\demo{{R}emark} The reader may verify by elementary calculus that
for $E \subset {\Bbb S}^{N-1}$ with
$$\hbox{area}(E) =(\frac12 - \varepsilon)\omega_N,$$
i.e., an $\varepsilon$ digression from the half-spherical cap,
we have
$$\frac12 - \varepsilon =\frac{1}{\sqrt{2\pi} }\int_h^\infty \exp{(-t^2/2)}
dt= \frac12 - \frac{1}{\sqrt{2\pi} }\int_0^h \exp{(-t^2/2)}dt\approx
\frac12 - h.$$
Hence $\varepsilon \approx h$ and therefore in Lemma 2.3
we will  have
$$\Lambda(h_+) + \Lambda(h_-) - 2 \geq 
C\Lambda''(0)(h_+^2 + h_-^2)\approx C(h_+^2 + h_-^2)\approx
C\varepsilon^2,$$
which gives (2.11).
\enddemo

\section{Proof of Theorem I}

First we need some definitions and notations.

\numbereddemo{Definition}
Set 
$$S_j(z,u)= \sup_{B(z,2^{-j})} |u|,$$
and define $\M (z,u)$ to be the maximal  subset of 
$\N$ (natural numbers) satisfying
 the following doubling condition 
$$ 4 S_{(j+1)}(z,u)\geq  S_{j}(z,u)\qquad \hbox{for all}\ 
\ j \in \M (z,u).\tag 3.1$$
\enddemo
Our aim is to prove that $S_j \leq C 2^{-2j}$, for all $j \in \M
(z,u)$
and some positive constant $C$.
An important observation at this point is that if $\M$ is empty 
 then we may
easily (by iteration) obtain the desired estimate.
Hence from now on we assume $\M \neq \emptyset$.

\proclaim{Lemma} Let $u \in P_1(z,M)${\rm .} Then there exists a constant $C_0=C_0(N)$ such that
$$S_j(z,u) \leq C_0 M 2^{-2j} \qquad \qquad \hbox{\rm for all}\  j \in \M (z,u). \tag 3.2$$
\endproclaim

\demo{Proof} By the remark following Definition 1.1 we may assume
that  $z$ is the origin.   Suppose the conclusion in the lemma fails. Then there exist
  $\{u_j \}$, $\{ k_j\}$ such that
$$S_{k_j}(0,u_j) \geq j 2^{-2k_j}, \tag 3.3$$
with $k_j \in \M(0,u_j)\neq \emptyset$.

 Now define $\tilde u_j$ as
$$\tilde u_j (x)= {u_j( 2^{-k_j}x)   \over S_{k_j +1}(0,u_j)  } 
\qquad \qquad \text{in } B(0,1).$$
Then $\tilde u_j$ satisfies the following properties:
$$ \V \Delta \tilde u_j \V_{\infty,B} \leq { 2^{-2k_j} \over S_{k_j + 1}(0,u_j) }
\leq {S_{k_j}(0,u_j) \over j S_{k_j + 1}(0,u_j)  }\leq { 4 \over j}\to 0,
\tag 3.4 $$
where the second inequality follows from (3.3), and the last inequality follows from
(3.1),
$$\sup_{B_{(1/2)}}|\tilde u_j|=1,  \tag 3.5$$
$$\Vert  \tilde u_j\Vert_{\infty,B} \leq {S_{k_j}(0,u_j) \over S_{k_j + 1}(0,u_j)}\leq 4,\tag 3.6$$
$$\tilde u_j(0)= |\nabla \tilde u_j|(0)=0  . \tag 3.7$$

Now by (3.4)--(3.7) we will have a subsequence of
$\tilde u_j$ converging in $C^{1,\alpha}(B)$ (see [GT]) to a 
nonzero harmonic function $u_0$, satisfying $u_0(0)=|\nabla u_0|(0)=0$.
For any fixed direction $e$ 
define  
$$v= D_e u_0, \qquad v_j =D_e u_j, \qquad \tilde v_j=D_e
\tilde  u_j.$$
Then,  for a subsequence,   $\tilde v_j$ converges
 in $C^{1,\alpha}(B)$  to $v$, where $\Delta v =0$.
Now according to Lemma 2.1 (since $v_j$ is harmonic in $\Omega_j$)
$${ 1 \over r^{2N}} \int_{B(0,r)} | \nabla v_j^+|^2 \int_{B(0,r)}
 | \nabla v_j^-|^2
\leq C \qquad \hbox{for all}\  r,j ,\tag 3.8$$
where $C$ depends on the $W^{2,2}$ norm of $u_j$ over the unit ball. By elliptic estimates this is uniformly bounded for
all $j$.  Making change of the variables in (3.8) and letting $r=2^{-k_j}$,
we will  obtain
$$ \int_{B(0,1)}| \nabla\tilde v_j^+|^2\int_{B(0,1)} | \nabla\tilde
v_j^-|^2\leq C\left( {2^{-2k_j} \over S_{k_j+1}}\right)^4
\leq  C\left( { 2^{-2k_j} \over S_{k_j}}\right)^4
 \qquad \hbox{for all}\ j ,$$
where in the last inequality we have used (3.1). Here, and in the
sequel,  $C$ is a  generic constant.
Next, invoking the Poincar\'e  inequality we may reduce (improve) the above to
$$ \int_{B(0,1)}| \tilde v_j^+ - M_j^+|^2\int_{B(0,1)} | \tilde
v_j^- - M_j^-|^2
\leq  C\left( { 2^{-2k_j} \over S_{k_j}}\right)^4
 \qquad \hbox{for all}\ j ,$$
where $M_j^\pm$ is the mean-value 
 for the functions $\tilde v_j^\pm$, on the unit ball.
 From here using   (3.3) 
we obtain, by  letting  $j$ tend to infinity, 
$$ 
\int_{B(0,1)} |  v^+ - M^+|^2 \int_{B(0,1)}  | v^- - M^-|^2 =0 , 
 \tag 3.9$$
where $M^\pm$ is the corresponding mean-value for $v^\pm$.
Obviously (3.9) implies that either of $v^\pm$ is constant.
Since also   $v(0)=0$, the constant must be zero.
In particular $v $ does not change sign.
But then 
 the maximum principle gives $v\equiv 0$, i.e., $D_e u_0 \equiv 0$.
Since $e$ is arbitrary we also have $u_0$ 
 is constant. 
Next using $u_0(0)=0$
we will  have $u_0\equiv 0$ which contradicts (3.5).
This proves the lemma. 
\enddemo

Next we will complete the chain in $\N$ for the estimate (3.2).

\proclaim{Lemma} Let   $C_0=C_0(N)$ be the constant in Lemma {\rm 3.2.}
Then 
$$S_j \leq 4 C_0 M2^{-2j} \qquad \hbox{\rm for all}\ j \in \N. $$
\endproclaim

\demo{Proof} Let $u \in P_1(z,M)$. Then obviously $S_1 \leq M$. Let
$j >1$ be the first positive integer such that 
the statement of the lemma  does not hold, i.e.,
$$S_j > 4C_0M2^{-2j}.\tag {3.10}$$
 Then
$$S_{j-1} \leq 4C_0M2^{-2(j-1)}=16C_0M2^{-2j}<
4S_j.$$
Hence $j-1 \in \M (u)$. By Lemma 3.2, then,
$$S_j\leq S_{j-1} \leq C_0M2^{-2(j-1)}=4C_0M2^{-2j},$$
which contradicts (3.10). The result follows.
\enddemo

From Lemma 3.3 we infer a uniform $C^{1,1}$ estimate for the class
$P_1(z,M)$.

\demo{Proof of Theorem {\rm I}} Let   $u$ be in $P_1(z,M)$ and set
$d(x)=\hbox{dist} (x, \partial \Omega)$. Then by Lemma 3.3
$$|u(x)| \leq C M d(x)^2 \qquad \qquad \hbox{for all} \ x \in B(0,1/2).\tag 3.11$$

Define now 
$$v(y) = {u(x + yd(x) ) \over d(x)^2} \qquad \qquad \hbox{ in } B(0,1).$$
 By (3.11)
$v$ is bounded on the unit ball, and by its definition it satisfies
$\Delta v =1$ in $B(0,1)$. Hence by elliptic estimates
$D_{ij} v(0) =D_{ij} u(x)$ is uniformly  bounded (independent of $x$).
This gives the result.
\enddemo

\section{Further auxiliary lemmas}

 In this and the next sections we will frequently  
 use the blow-up of functions;
i.e.\  for a  given $u$ we consider $u_r(x)=u(rx)/r^2$ and
let $r $ tend to zero, through some subsequence. It is, however, not
clear whether the blow-up (the limit function) will not be the zero
function.
Indeed, if $u(x)=o(|x|^2)$,
 then any blow-up will be identically zero. To prevent this we need a 
nondegeneracy from below,  asserted in the following remark.
\demo{{R}emark} The function $u$ in (1.1) satisfies
$$\sup_{B(x^0,r)} u \geq u(x^0) + C_Nr^2, \tag{4.1}$$
for all $x^0 \in \overline \Omega$. Here  $C_N=1/2N$,
if $u(x^0) \geq 0$, and $C_N$ is somewhat smaller  if $u(x^0)<0$.
Also  $r$ is small enough so that $B(x^0,r) \subset B(0,1)$.
\enddemo

The proof of this is given  in [Ca2] in  the case $u > 0$ or $x^0 \in
\partial \Omega$. The
general case is proven as follows.
Suppose $u(x^0) \leq 0$ and $x^0 \in \Omega$.
 If there is $x^1 \in B(x^0,r/2)\cap \partial
\Omega$,
then we apply the above to $u$ in $B(x^1,r/2)$ to obtain
the same estimate with $C_N=1/8N$. So suppose
 the set $B(x^0,r/2)$ contains no
free boundaries. 
we can apply the mean value theorem for harmonic functions to 
$u(x) - |x- x^0|^2/2N$ to obtain
$$\int_{B(x^0,r/2)}( u  -u(x^0)) \ dx= cr^{2+ N},$$
where $c$  depends on $N$ only.
From here one obtains (4.1).

Now,  by (4.1) and Theorem I we may assume that any blow up 
of functions in $P_1(0,M)$ remain
in  the class.
Our next definition will be used frequently here and later in Section 6.

\demo{$\varepsilon$\/{\rm -}\/close}  We say two functions $f$ and $g$ are
$\varepsilon$-close to each other in a domain $D$ if 
$$\sup_{x \in D} |f(x) -g(x)| < \varepsilon.$$
\enddemo

\demo{Blow\/{\rm -}\/up limit} A blow-up limit $u_0$ is a uniform limit on
compact subsets of $\Bbb R^N$
$$u_0(x) =\lim_{j \to \infty} \frac{u(r_j x)}{r_j^2}$$
where $u \in P_1(0,M)$ and $r_j \to 0$.
The function $u$ may even change for different~$j$.
\enddemo

\demo{Flat points}
We say   $\partial \Omega$  is flat at the origin  if there is a blow
up $u_0$   such that the set where  $u_0 = 0$
is a half-space. 
\enddemo

\demo{Half\/{\rm -}\/space solutions}
A half-space solution  $u_0$ is a global solution that has the representation
$(\max(x_1,0))^2/2$ in some coordinate system.
\enddemo

For the readers convenience  and for   future reference 
we will  recall and explain some general (known) facts.
These will be crucial in the rest of the paper. 
We recommend that a  reader unfamiliar  with  such problems carefully
verify these facts.

\demo{General Remarks}  

{a)}  
 By (4.1) and the techniques of [Ca4, Lemma 6] one may show that
the set 
$$\Omega \cap \{ x: \ |\nabla u(x)|<\varepsilon \ \}$$ has
volume less than $C\varepsilon$, with universal $C$.
 One uses only the $C^{1,1}$ property of the
solution, and not the nonnegativity of the function. 

This implies, in particular, that 
 an $\varepsilon$-neighborhood
$$K_\varepsilon=\{x: \ \hbox{dist}(x,\partial \Omega)< \varepsilon\}$$ 
of the free 
 boundary has Lebesgue measure less than $C\varepsilon$,
i.e.\ volume$(K_\varepsilon) <C \varepsilon$, with universal $C$.

Using covering arguments, we conclude that 
the  free boundary $\partial \Omega$
has locally finite $(N-1)$-Hausdorff measure; see [Ca4] for details.

\medbreak
{b)}  Let $u_j$ be a blow-up of $u$ and
suppose $u_j$ converges to $u_0$ in $C_{\rm loc}^{1,\alpha}(\R)$.
It follows that
 $$\{u_0=|\nabla u_0|=0\} \supset \varlimsup  \{u_j=|\nabla u_j|=0\},$$
(see e.g. ([Ca2], [KS1]));   $\varlimsup$ denotes the limit set of
all sequences $\{x_j\} , x_j \in \{u_j=|\nabla u_j|=0\}$.

Next let $u_0$ be a blow-up of a sequence $u_j$ 
and suppose  $u_0=0$  in  $B(x^0,r_0)$ for some $x^0$, and $r_0$.
Uniform convergence and (4.1) imply that
  $u_j=0$ in $B(x^0,r_0/2)$ 
for large $j$.

A consequence of this  is  the following
$$\hbox{interior}
(\{u_0=|\nabla u_0|=0\}) \subset \varlimsup  \{u_j=|\nabla u_j|=0\}.$$

\medbreak
{c)} From a) and b) above we infer an $L^p$ convergence of
the second derivatives of $u_j$ to $u_0$; i.e.
$$D_{ik} u_j  \to D_{ik} u_0 \qquad \hbox{in } L^p_{\rm loc}-\hbox{norm},$$
for $1<p<\infty $. This depends on the
fact that $\Delta u_j$ and $\Delta u_0$ differ only 
inside an $\varepsilon$-neighborhood of the free boundary $\partial
\Omega_0$;
i.e.\ on a set of Lebesgue measure~$\varepsilon$.

\medbreak
{d)}
From a) and b) we may also deduce that if $u\in P_1(0,M)$ is $\varepsilon$-close to 
 a half-space solution $h=(\max(x_1,0))^2/2$, say, then
 $u\equiv 0$ in $$B(0,1/2) \cap \{ x_1 < -C\sqrt{\varepsilon}\},$$ 
for some constant
$C>0$. 
We sketch some details.  Let us suppose 
$x^0 \in  B(0,1/2) \cap \{x_1 < 0\} \cap \Omega$.
Choose $r=|x_1^0|$.
Then by (4.1) and the closeness of $u $ to the half-space solution $h$
we have
$$2\varepsilon \geq C_Nr^2.$$
Observe that $|u|<\varepsilon$ in $B(x^0,r) \subset \{x_1 < 0\}$.
Hence if $x^0_1<-\sqrt{2\varepsilon/C_N}$ then $x^0$ cannot be a point of
$\overline  \Omega$.
This implies that   $u \equiv 0$ on  the set
$\{x_1<-\sqrt{2\varepsilon/C_N}\}$.

\medbreak {e)}
A consequence of this is the following simple fact:
Let $u \in P_1(0,M)$  and suppose  the origin  is a flat point 
with respect to some blow-up sequence $u_{r_j}$. Then in $B(0,1/2)$,
$u_{r_j}$ is $\varepsilon$-close to a half-space solution for small enough~$r_j$.
\enddemo
  
\proclaimtitle{essentially due to Spruck [Spk]}
\proclaim{Lemma}   Let $ u \in P_1(0,M)$.
 Then any blow\/{\rm -}\/up $u_0$ of $u$ is a  homogeneous function of degree two{\rm ,}
 and the set $\{u_0 =|\nabla u_0|= 0\}$ is a cone{\rm .}
\endproclaim

The proof of Lemma 4.1 in the global case ($P_\infty(0,M)$) is given in details in [KS2].
The proof of the local case is similar  and therefore omitted.

\proclaim{Lemma}  Let $ u \in P_1(0,M)$ and suppose 
$CD_1 u - u \geq - \varepsilon_0$ in $B(0,1)$  for some $\varepsilon_0 >0$
and $C>0${\rm .}
Then $CD_1 u - u \geq 0$ in $B(0,1/2)${\rm ,} provided $\varepsilon_0$ is small
enough{\rm .}
In particular{\rm ,}  if $u$ is close to the   half\/{\rm -}\/space solution
$\max(x_1,0)^2/2$ in $B(0,1)${\rm ,} then {\rm (}\/by integration and {\rm d)} in General
Remarks\/{\rm )}  $u \geq 0$ in $B(0,1/2)${\rm .}
\endproclaim

\demo{Proof}   Suppose the conclusion of the 
lemma fails. Then there is a  $u \in P_1(0,M)$ with
$$CD_1 u(x^0) - u(x^0) < 0 ,\tag{4.2}$$ 
 for some $ x^0 \in   B(0,1/2)$.
Let
$$w(x) = CD_1 u(x) - u(x) + \frac{1}{2N}|x - x^0|^2.$$
Then $w$ is harmonic in $\Omega \cap B(x^0,1/2)$, $w(x^0)<0$ (by (4.2)) and
$w \geq 0$ on $ \partial \Omega$. Hence by the maximum principle
the negative infimum of $w$ is attained on $\partial B(x^0,1/2)$.
We thus obtain
$$-\varepsilon_0 \leq \inf_{\partial B(x^0,1/2) \cap \Omega}
(CD_1 u - u ) \leq -{1\over 8N},$$
which is a contradiction as soon as $\varepsilon_0 < 1/(8N)$.
The second part, that $u \geq 0$ in $B(0,1/2)$ for $u$
near to half-space solution, 
 follows by d) in General Remarks and
integration. The reader should notice that
$\varepsilon$-closeness
for $u$ to a half-space solution implies, by elliptic estimates,
 that the gradient of $u$
also becomes $(C\sqrt{\varepsilon})$-close to the gradient of the half-space
solution. Here $C$ is a universal constant. We leave the details to
the reader.  This proves the lemma.
\enddemo

In the next lemma we will apply the following  simple fact, which we
  formulate as a remark.

\demo{{R}emark} 
  Suppose  $u_0$ is a blow-up
solution of  $u$, and  $\delta_1(u_0) >0$. Then $u_0$ is a degree two
 homogeneous global solution and the 
interior of $ \R \setminus \Omega$ is a   nonvoid cone. 

Indeed  by Lemma 4.1,  $u_0$ has the   mentioned properties and
  $\R \setminus \Omega$ is a cone. If  the cone has an
empty interior, then by General Remark a), ${\Bbb R}^N \setminus \Omega$
has Lebesgue measure zero. Therefore
by Liouville's theorem 
 $u_0$ is a polynomial of
degree two and hence for all $r$ we will  have
 $\delta_r(u_0)=0$, which is a  contradiction.
\enddemo

\proclaim{Lemma} Let $u \in P_1(0,M)${\rm ,} and suppose
$\varlimsup_{r \to 0} \delta_r (u) > 0${\rm .} Then $u \geq 0$ in some neighborhood of
the origin{\rm .}
\endproclaim
\demo{Proof} Let $\{r_j\}$ be a decreasing sequence such that
$$\varlimsup_{j \to \infty} \delta_{r_j}(u) > 0.$$
We blow up the function $u$ through $r_j$ to obtain, by Lemma 4.1,
 a global homogeneous
solution $u_0$ of degree 2 in  $\R$, with 
$$\delta_1(u_0) >0. $$
Hence  by the discussion preceding this
lemma, $\R \setminus \overline \Omega_0$ is a nonvoid cone.

 Now choose a direction $e$  and consider the
monotonicity formula for $D_e u_0$. By degree two homogeneity of
$u_0$ we have  $|\nabla D_e u_0|$ is homogeneous of degree zero. By
scaling,  this implies that $\varphi(r, D_e u_0)$ must
 be a positive constant  unless
one of the functions $(D_e u_0)^\pm$ is zero.
 To exclude the first case, observe that by Lemmas 2.2--2.3 if 
$\varphi \neq 0$ then it is strictly
monotone  since (by the condition $\delta_1(u_0) >0$) at least one of the sets $\Gamma_\pm$ cannot be a  half-sphere.
 In the second case we will 
have  $D_e u_0 \geq 0$  (or $\leq 0$) 
for any fixed  direction $e$.

Let $x^0$ be a fixed  point in $\Omega_0$ and suppose $|\nabla
u_0(x^0)|\neq 0$. Set 
$$\nu =\frac{\nabla u_0(x^0)}{|\nabla u_0(x^0)|}.$$
Then for any directional vector $e$ orthogonal 
to $\nu$ we will have  $D_e u_0(x^0)=0$. 
Moreover, $D_eu_0 \geq 0$ in $\Omega_0$ (or $\leq 0$). It is harmonic
there and it takes a local minimum (or maximum).
 Hence by the minimum (or maximum) principle $D_e u_0 =0$
in $\Omega_0$. This implies that $u_0$ is independent of the directions~$e$ orthogonal to $\nu$ and depends only on
the direction $\nu$. Therefore
$u_0$ is one-dimensional in each connected component.
Since the only one-dimensional solutions are half-space solutions
(in each connected components) there must be at most two
connected components  which are half-spaces.  
Now  the assumption  $\delta_1(u_0)>0$ implies
that there must be at most one connected components  which is a
half-space.

In particular, near the origin, $u$ is close to a half-space solution.
Therefore the origin is a flat point.
Hence  Lemma 4.2  and part e) in General Remarks 
 give  the result.
 This proves the lemma.  
\enddemo  

Our next result will not be used in this paper, however, we include
this for future references.
First we notify the reader of the following obvious fact.
\demo{{R}emark}
 Suppose  $u \in P_1(0,M)$ and the origin is a 
point of zero upper Lebesgue density for the complement of $\Omega$.
 Then any blow-up of $u$ at the origin is a polynomial of degree two.
\enddemo

\proclaim{Lemma} Let $u \in P_1(0,M)$  satisfy $\varlimsup_{r \to 0}
\delta_r(u) =0$  and suppose for some blow\/{\rm -}\/up sequence with limit
$u_0${\rm ,} $D_e u_0 \neq 0${\rm .}
Then  there is a constant $C_e$ such that 
  $\varphi (r,D_e u) \geq C_e >0$  for all $r<1${\rm .}
\endproclaim

\demo{Proof} Blow up $u$ to obtain a global solution $u_0$, in $\R$ with 
${\rm int}(\Omega^c_0)=~\emptyset$ (see the discussion preceding this lemma).
 Hence $u_0=P$, a degree two polynomial.
 For any direction $e$ nonparallel to the kernel of $P$ we have
$D_e u_0 =D_e P$ is a nonconstant  linear function. Next
$$\varphi(r) \geq {1 \over r^{2N}}\int_{B_r} |\nabla(D_e u)^-|^2\int_{B_r} |\nabla(D_e u)^+|^2.
$$
Scaling and letting $r $ tend to zero we will  have (see c) in General Remarks)
$$\varphi(0) \geq \int_{B_1} |\nabla( D_e P)^-|^2\int_{B_1} |\nabla (D_e P)^+|^2 \geq
C_e ,\tag{4.5}$$
where in the last inequality we have used the fact that $D_eP$ is a
nonconstant linear polynomial. Now (4.5) together with the
monotonicity formula
gives the result. 
\enddemo

\demo{{R}emark} Lemma 4.4 can be stated in a more accurate way. We
can, with suitable choice of $e$, make the constant $C_e$ uniformly bounded from below for
the
whole class $P_1(0,M)$. Here is how:
Rearrange the coordinate system such that
the polynomial $P$ in the proof of Lemma 4.4 has the representation
$$P=\sum_{i=1}^m a_i x_i^2 \qquad \sum_{i=1}^m a_i =\frac12,$$
where $m \leq N$.
Let now $e_j$ be the standard coordinate system and set
$$e=\frac{(e_1 + \cdots +  e_m)}{\sqrt{m}}.$$
Then 
$$D_eP=\frac{2}{\sqrt{m}} \sum_{i=1}^m a_i x_i,$$
and
$$|\nabla D_e P|^2 =\frac4m \sum_{i=1}^m a_i^2
\geq \frac{4}{m^2} (\sum_{i=1}^m a_i)^2 =\frac{1}{m^2}\geq \frac{1}{N^2}.$$
Hence for $e$ as above 
$$\varphi(0,D_e u) \geq \left(\frac{\hbox{vol}(B_1)}{2N^2}\right)^2.$$
\enddemo

\section{Proof of Theorem II}

 We remark that in the definition of the global solutions
$P_\infty(0,M)$, we require that the functions have quadratic growth with uniform constant $M$. It is noteworthy that this restriction is not 
superfluous  as
there are examples of  solutions to (1.1) in the entire space  with 
$\R\setminus \overline \Omega$ nonvoid and with $u$ of polynomial or even exponential growth;
see [Shap1].

However, we will use Theorem II  in connection with blow-up functions,
 i.e.,
we consider blow-up of functions in $P_1(0,M)$ and these, by Theorem I,
have quadratic growth near the origin.   Hence any blow-up of such
functions will also be of quadratic growth in the entire space.

The reader should observe that, by the assumption $\delta_r(u) >0$
for some $r>0$ in Theorem II, and the discussion in the remark
preceding Lemma 4.2,
${\Bbb R}^N \setminus \Omega$ has nonempty interior. It suffices to
show that $u \geq 0$, since  by [Ca2] it follows that $D_{ee}u
\geq 0$ on $\partial \Omega$  and by scaling we also have this
property in $\Omega$ (see
the details of this scaling argument in [KS2]). We will prove

\specialnumber{II'}
\proclaim{Theorem}  If $u \in P_{\infty}(0,M)$ and ${\Bbb R}^N
\setminus
\Omega$ has nonempty interior{\rm ,} then $u \geq 0$.
\endproclaim

We split the proof in three cases.
\medbreak
 
{\it Case} 1. $\R \setminus \Omega$ is bounded.

 \smallbreak
{\it Case} 2.  $\varlimsup_{R\to \infty }\delta_R (u) >0$.
\smallbreak
 
{\it Case} 3. $\varlimsup_{R\to \infty }\delta_R (u) =0$.
\pagebreak

\demo{{R}emark} The reader should observe that, by the assumption 
$\delta_r (u) >0$ (for some $r>0$) in Theorem II, and 
the discussion in the remark preceding Lemma 4.2,
 $\R \setminus \Omega$   has nonempty  interior. 
 The second part of  Theorem II, the fact that 
$\varlimsup \delta_R (u)>0$  implies
$u$ is a half-space solution, is included in the proof of Case 2.
\enddemo
\demo{ Proof of Case {\rm 1}} Let $U$ be the Newtonian potential
of the complement of $\Omega$  with constant density  and
such that $\Delta U =\chi_{\Omega^c}$,
i.e.\ ($U=c| .  |^{2-N} \ast \chi_{ \Omega^c}$;
for $N=2$ we take the logarithmic kernel).
Then, since $\partial \Omega$ has zero Lebesgue measure,
 $\Delta( U + u ) = 1$ almost everywhere in $\R$, and it has quadratic growth.
 Hence by Liouville's theorem 
it is a  second degree polynomial $P$.

By translation and rotation  we may assume $P=\sum a_j x_j^2 + d$;
observe that, by the imposed rigid motion, 
the origin is not necessarily on $\partial \Omega$ any more.
 From here we will have 
$u = P -U$ in $\Omega$.
 Let now $v = x\cdot \nabla u - 2u$.
Then by homogeneity of $P-d$ we will have 
$$v= -2d - x\cdot \nabla U + 2U \quad \to\quad  
C 
\qquad \qquad \hbox{ as } x \to \infty,$$
 and $v=0$ on $\partial \Omega$. 
If $N=2$ then $C=+\infty$ and if $N \geq 3$ then $C=-2d$.

Suppose first $C<0$. Then by the
maximum principle $v<0$. Now  by elementary calculus we will imply
 $u(rx)/r^2$ is decreasing in $r$.
On the unbounded
cone-like set
$K=\{rx: \ x \in \R \setminus \Omega,\  r \geq 1
\ \}$,
we will have
 $u \leq 0$. By subharmonicity, this implies
$$0=u(x^0)\leq \hbox{vol}(B) \int_B u,
\qquad \qquad \hbox{for all } \ x^0 \in \partial \Omega \cap \hbox{int}(K),$$ 
 where $B$ is
a ball in $K$ with center $x^0$. 

Since also $u\leq 0$ in $K$,
we will have  $u\equiv 0$ in $B$,
 i.e.\ $\partial \Omega \cap \hbox{int}( K)= \emptyset$.
Therefore $K \subset \R \setminus \Omega$. But then $\R \setminus \Omega$
is unbounded, which is a contradiction. Therefore $C>0$.

Now   $C>0$ implies, by the maximum principle, that $v >0$,
i.e., $u(rx)/r^2$ is increasing in $r$.
 Then $u\leq 0$ on the truncated cone
 $$K=\{rx: \ x \in \R \setminus \Omega,\  r \leq 1
 \}. $$   A similar argument as above then shows that $u=0$ on $K$.
Hence  we conclude  that the set  $\R \setminus \Omega$ has
positive Lebesgue density at the origin.

 Now
two cases arise: (a) the origin is in the interior of $\R \setminus \Omega$,
and (b) the origin is on  $\partial \Omega$.
The first case, along with the monotonicity of $u(rx)/r^2$
will result in the positivity of $u$. As to the second case,
we  observe that  since $K \subset \{u=0\}$, 
$\varlimsup_{r \to 0} \delta_r (u) >0$ and by Lemma 4.3, 
$u \geq 0$ near the origin.
Again by the monotonicity of $u(rx)/r^2$ in $r$ we will have  $u \geq 0$
in $\R$.   \hfill\qed
\enddemo

A simple alternate, but indirect, proof for   Case 1  goes 
as follows: Since the Newtonian potential $U$
is a polynomial in $\Omega^c$, it follows from  [DF]
 that $\Omega$ is the exterior of an ellipsoid
and from [Shah] it follows that $u \geq 0$. 

\demo{ Proof of Case {\rm 2}} We first blow up $u$ through a  sequence
$\{R_j\}$
($R_j \to \infty$) such that $\lim\delta_{R_j} (u) \neq 0$ to obtain a
global homogeneous solution $u_0$ of degree two (see [KS2,  Lemma 2.5]).
A similar argument as that in the proof of Lemma 4.3 then  implies
that $u_0$ is a half-space solution. In particular  this  implies
$$\varphi(\infty,D_eu):=\lim_{R_j \to \infty}\varphi(R_j,D_e u)=
\lim_{R_j \to \infty}\varphi(1,D_e u_{R_j})=\varphi(1,D_e u_0)=0,$$
 for all directions
$e$. 
 The last equality depends on the fact that  a half-space
solution $u_0$ is monotone in every  direction, and thus
either $(D_eu_0)^+\equiv 0$ or $(D_eu_0)^-\equiv 0$.
Observe also that  we have used
 the convergence in $W^{2,2}$; see c) in General Remarks.

Now  by the monotonicity
formula, 
$$\varphi(r,D_e u)\leq \varphi(\infty,D_e u)= 0,$$
 for all vectors  $e$. Hence $D_e u \geq 0$ (or the reverse); i.e.
 $u$ is monotone in
all directions.  As in Lemma 4.3 it follows that $\nabla u$ is
parallel
at any two points of a component of $\Omega$, and hence that $u$ is a
half-space solution. \hfill\qed
\enddemo

\demo{ Proof of Case {\rm 3}} Since $\R \setminus \Omega$ is unbounded and 
$\varlimsup_{R \to \infty} \delta_R (u) =0$ we may conclude that there
is a  blow-up $u_j(x)=u(R_jx)/R_j^2$
at infinity ($R_j \to \infty$) with  a subsequence converging to a polynomial $P$ in $\R$, and
that $P$ is independent of some of the variables. 
Indeed, by the assumptions in this case
there is an unbounded  sequence $x^j \in \partial \Omega$. 
Therefore  we may take
$R_j=|x^j|$ and obtain that $P$ vanishes, along with its gradient,
at the origin and at  some other point on  the unit sphere. Hence by 
homogeneity (Lemma 4.1) the same is true 
on the whole line generated by these points, and thus the above conclusion.

Suppose $D_1 P=0$. Then for any point $x^0 \in \partial \Omega$ we may consider the monotonicity formula
$\varphi(r,D_1u,x^0)$, which by Lemma 2.1 is nondecreasing in~$r$.
Hence for $r$ fixed we choose $R_j \geq r$ to obtain
$$\varphi(r, D_1 u,x^0) \leq \varphi(R_j, D_1 u,x^0) = \varphi(1, D_1
u_j,x^0)  \to \varphi(1, D_1 P,x^0) =0,$$
as $R_j \to \infty$. Here  we have used that $D_1 u_j \to D_1 P =0$.
This will imply that  $\varphi \equiv 0$, i.e.\ either of $(D_1 u)^+$ or $(D_1 u)^-$ is zero.
Suppose  $D_1 u \geq 0$ (the other case is treated similarly).
From here we want to deduce that $D_{11} u \geq 0$.
To do this we set
$$-C:=\inf_{\Omega} D_{11} u ,$$
which is bounded because $u$ is in $P_\infty(0,M)$.

Let $x^j$ be a minimizing sequence, i.e.,
$$-C=\lim_j D_{11} u(x^j).$$
We consider a blow up at $x^j$ and with
$d_j=\hbox{dist}(x^j,\partial \Omega)$. Hence we define
$$u_j(x)=\frac{u(d_jx + x^j)}{d_j^2} \qquad \hbox{in } B(0,1).$$

Now by compactness, for a subsequence
we will  have
$ u_j \to u_0,$
where $u_0$ is a global solution.
Also 
$D_{11} u_j \to D_{11}u_0$ uniformly in $B(0,1/2)$.
The latter depends on the fact that $\Delta u_j =1$ in
$B(0,1)$. Hence  we will  have a global solution $u_0$
with $ \Delta u_0 =\chi_{\Omega_0}$
where $\Omega_0 = \Omega(u_0)$.
From the minimal properties of the sequence $x^j$ we also deduce
that for all $  x \in B(0,1/2),$
$$D_{11}u_0(0) = -C \leq
\lim_j D_{11} u_j (x) = \lim_j D_{11} u (d_jx+x^j)= D_{11} u_0 (x). $$
 Hence by the maximum principle $D_{11}u_0 \equiv -C$ in the connected
component of $\Omega_0$ that contains the unit ball, we call this $\Omega'$.
Observe also that the free  boundary $\partial \Omega_0$ is
nonempty. Therefore we assume that $z=(z_1, \cdots , z_N)\in \partial
\Omega_0$.
Integration 
 gives that in $\Omega'_0$ we have the following representation
$$D_1u_0(x) = -Cx_1 + g(x_2, \cdots , x_N).$$
Since also 
$D_1 u_0 = \lim_j D_1u_j \geq 0 $ we will  have
$x_1 \leq g(x')/C$ (of course we assume that $C >0$ otherwise there is
nothing to prove).
 In particular this means that any ray $l_x$ emanating at $x \in
 \Omega'$ and parallel to the $x_1$-axis hits $\partial \Omega'$.
Now the component   $\Omega' $ of $\Omega_0$
is  the under-graph of  the function $x_1=g/C$  in the
$x_1$-direction.  $D_1u_0 \geq 0 $ implies that
in the negative $x_1$ direction (inwards to $\Omega'$)
$u_0$ is decreasing. Since it is also  zero on $\partial \Omega'$
 it becomes  nonpositve  in $\Omega'$. This contradicts the
 nondegeneracy (4.1).

 Therefore we conclude $C\geq 0$ and hence
$$D_{11}u \geq 0 \qquad \hbox{in } \Omega.$$

This along with $D_1 u \geq 0$ 
implies that $u \geq 0$ on lines  which hit $\R \setminus \Omega$
 and are
 parallel to the $x_1$-axis.
Our goal will be to prove that for any $x$
in $\R$,
$$\lim_{ m \to \infty }  u(x_1 - m, x_2 , \cdots , x_N)\geq 0,$$
 which together with $D_1 u \geq 0$ implies  $u(x) \geq 0$.
Now let 
$$u_m(x)= u(x_1 - m, x_2, \cdots, x_N)\qquad  (m=1,2, \cdots )$$
be a family of translations of $u$. Since  the  negative $x_1$-axis is
in $\R \setminus \Omega$ (observe that $D_1 u \geq$ and $u(0)=0$)
 and $|u(x)|\leq M (|x|^2+1)$ we  deduce
that $|u_m (x) | \leq C_0 R^2$ on $B(0,R)$ ($C_0$ independent of $m$).
This in particular implies that $u_m$ is a bounded sequence;
hence there is a converging subsequence with a limit function $u_0$
in $\R$. It is then elementary to see
that the  function $u_0$ is also a global solution. 

Next using
the nonnegativity of $D_{11}u$ and $D_1 u$ we infer that
the nonnegative monotone function $D_1 u$ has a limit at
$-\infty$; i.e.,

$$D_1 u_0 (x) = \lim_{ m \to \infty } D_1 u(x_1 - m, x_2 , \cdots , x_N)= \hbox{ constant}=0,$$
\vglue6pt\noindent 
and hence  $u_0$ is cylindrical (($N-1$)-dimensional) in
 the $x_1$-direction. If $u_0 \equiv~0$ there is nothing to prove;
therefore we assume that $u_0$ is nontrivial  and
one lower-dimensional. Also the set $\{u_0 \equiv 0\}$ has nonempty
interior since it contains the projection of the interior of the set
$\{u\equiv 0\}$.
But   one of the cases may occur for this new lower-dimensional
function  and we repeat the argument (if the third
case occurs) until we obtain a one-dimensional problem. 

Now  the reader may
verify, using elementary calculus, that the one-di\-mensional global solution
is nonnegative. This proves  $u \geq 0$. \hfill\qed
\enddemo
\vglue4pt

\section{Proof of Theorem III} 

\vglue4pt

In this section we will give a proof of Theorem III. The proof 
is based on two lemmas that are consequences of already established
results,  namely
 Theorem II,  the techniques in Lemma 4.2, and the main result in [Ca2].
We first use Theorem II in conjunction with [Ca2] to
obtain uniform  flatness for the free boundary  of  $u \in P_\infty (0,M)$
provided $\delta_1 (u) \geq \varepsilon$. Next we show that for 
$u \in P_1(0,M)$, there is a uniform  neighborhood of $0$ such that
$u \geq 0$ in that neighborhood provided  $u$  is
 close  to a half-space solution ($\varepsilon$ close as in Section~4).

The first author's original result for nonnegative solutions
in the class $P_1(0,M)$  provides us with the following lemma;
see [Ca2] or [Ca4].

\phantom{bummer}

\proclaim{Lemma} Given a  positive number $\varepsilon${\rm ,} there
 exists $t_\varepsilon$  such that
if\break  $u \in P_\infty(0,M)$ and $ \delta_1(u) \geq \varepsilon${\rm ,}
 then  in $B(0,t_\varepsilon)$ the boundary of $\Omega$
is the graph of a $C^1$ function {\rm (}\/uniformly for the class\/{\rm )}
and $u$ is $(\varepsilon t^2_\varepsilon)$\/{\rm -}\/close to a half\/{\rm -}\/space solution there{\rm .}
\endproclaim
\vglue6pt

Now a simple proof of  Lemma 6.1 can be given based on compactness,
 a contradictory  argument, and the main result in [Ca2].
 
The proof of the next lemma
follows   the same lines as that of the proof of Lemma 4.2. See
also General Remarks. We omit the proof.

\proclaim{Lemma} Let $\varepsilon ,  s >0${\rm ,} and suppose 
$u \in P_1(0,M)$ is $(\varepsilon s^2)$\/{\rm -}\/close to a half\/{\rm -}\/space
solution $(\max(x\cdot e,0))^2/2$ in $B(0,s)${\rm .} Then  in $B(0,s/2)$
we have
$$s D_eu - u \geq 0,  \tag{6.1}$$
provided $\varepsilon$ is small enough{\rm .}
 In particular{\rm ,} by integration and {\rm d)} in General Remarks{\rm ,}
 we have $u\geq 0$ in $B(0,s/4)$ 
{\rm (}\/again if $\varepsilon$ is small enough\/{\rm ).}
\endproclaim

\demo{Proof of Theorem {\rm III}}  First we claim  that for given
 $\varepsilon >0$ there exists
$$0<r_\varepsilon<t_\varepsilon, \qquad 
\hbox{(where $t_\varepsilon$ is as in Lemma 6.2)}$$
such that if for $u \in P_1(0,M)$ we have
$\delta_{r_0}(u)\geq \varepsilon $ for some $r_0<r_\varepsilon$ 
then  $u$ is $(2\varepsilon r_0^2t^2_\varepsilon)$-close to a half-space
solution $(\max(x\cdot e,0))^2/2$
in $B(0,r_0  t_\varepsilon  /2)$.
 Suppose for the moment that  this is true.
By Lemma 6.2 it follows 
that $u \geq 0$ in $B(0,r_0^2/4)$. Moreover
$(2\varepsilon r_0^2t^2_\varepsilon)$-closeness 
to a half-space solution also implies (see d) in General Remarks) that
$u\equiv 0$ in $B(0,r_0^2/4)\cap \{x_1 < -C\sqrt{\varepsilon}r^2_0\}$.
Hence for small $\varepsilon$,  [Ca2] or [Ca4] 
applies to conclude that
in $B(0,c_0r_0^2)$,    $\partial \Omega$
is the graph of a $C^1$-function  with a uniform $C^1$-norm.
Here $c_0$ depends only on $M$ and $N$.

 The modulus of continuity $\sigma(r)$ is then defined
 by the inverse of the  relation 
$\varepsilon \to r_\varepsilon$.

Now to complete the proof we need to
show the $(2\varepsilon r_0^2t^2_\varepsilon)$-closeness of $u$ to a
half-space solution in $B(0,r_0t_\varepsilon/2)$.
Suppose this  fails. Then  for each $r_j \searrow 0$ 
there exists $u_j \in P_1(0,M)$
with $\delta_{r_j}(u_j)\geq \varepsilon $  such that 

$$\sup_{B(0,r_jt_\varepsilon/2)}|u_j - h| > 2\varepsilon r_j^2t_\varepsilon^2,$$
\vglue6pt\noindent 
for all half-space solutions $h$.

Set $v_j=u_j(r_jx)/r_j^2$. Then $v_j \in P_{1/r_j}(0,M)$,
$\delta_1(v_j) \geq \varepsilon$, and 

$$\sup_{B(0,t_\varepsilon/2)}|v_j - h| > 2\varepsilon t_\varepsilon^2,$$
\vglue6pt\noindent 
for all half-space solutions $h$.

 Now,  for a subsequence and in an 
appropriate space, $v_j$ converges to a  global solution $v_0$ which
 on one side satisfies the hypotheses of Lemma 6.1,
  but on the other
side 
$$\sup_{B(0,t_\varepsilon/2)}|v_0 - h| \geq 2 \varepsilon t_\varepsilon^2,$$
\vglue6pt\noindent 
for all half-space solutions $h$.
 By Lemma 6.1,
this is a contradiction. Hence the claim holds.

This completes the proof of the theorem.
\enddemo

\section{Application to the Pompeiu type  problem}

Theorems I and  III    above 
 apply to more general operators of the form
$$L(u) = \sum D_i ( a_{ij}D_j u) + a(x) u,$$
with $a_{ij}$ and $a$ in the $C^\alpha$-class and
$a_{ij}(0)=\delta_{ij}$ (Kronecker's delta function). 
We refer to  [Ca3] for these more general operators, and
with $a(x)\equiv 0$. Here we treat the case of $\Delta  + a(x)$ with
$a(x)$ constant.

To show this we need  to make sure that the monotonicity formula works
for these operators. However, this is not true in general, as one may
easily observe. On the other hand, the full strength of
the  monotonicity formula is not used
in our analysis.

 The monotonicity formula is used  either  locally near the origin or
globally when the blow-up functions are considered.
For the first case we simply will show that the  operator $L$ will admit
a monotonicity formula which is almost increasing. 

 The second
case is even simpler, as the blow-up of any function in the class
$P_1(0,M,L)$ ($L$ denotes the dependent on the operator) will converge
to $P_\infty (0,M,\Delta)$, i.e., after rotation the blow-up will satisfy
$\Delta u_0 =\chi_{\Omega_0}$  and we are back to the  same situation
as before.

 As to the technique in Lemmas 4.2 (and  6.2) we observe that
for the case $\Delta u + u =1 $  we will have
$\Delta w = -Cr D_e u - 1 + u + c_0 \leq 0$ if $c_0$ is small
enough and the lemmas work. 
Spruck's
theorem also works with almost no changes. We leave the details to the
reader.

 A particular problem, where the Laplacian is replaced with
the Helmholtz operator,  is the classical and well-studied 
Pompeiu problem which can be stated  as follows.
Suppose there exists a bounded  domain $\Omega \subset \R$ and a function
$u$ satisfying
$$\Delta u + u =\chi_\Omega\quad
\hbox{in } \R, \qquad  u=0 \quad\hbox{in } \R \setminus \Omega.\tag {7.1} $$
 Does it follow that
$\Omega$ is a ball?  This question, put forward 
somewhat differently by  Dimitrie Pompeiu,
  has puzzled many mathematicians during
the past 50 years. There are  many partial results for this problem,
which we will not discuss. For  background  in the Pompeiu problem and
 for further references and how the formulation of Pompeiu is related
 to  the one here see  [W].

 From our   perspective it is interesting to know whether the boundary
of such domains are analytic. Williams [W] proved that
if the boundary of $\Omega$ in (7.1) is Lipschitz then it is analytic.

It does follow from our results, with the  modification  below,
that under the thickness condition of Theorem III, the boundary
of such domains are analytic.
We formulate this in the following theorem.
\specialnumber{IV}
\proclaim{Theorem} 
Under the thickness assumptions of Theorem {\rm III,}
any domain $\Omega$ which admits a solution $u$   to {\rm (7.1)} has an
analytic boundary{\rm .}
\endproclaim  

Now to show the almost monotonicity of the function $\varphi$
we need to define the following functions,
$$I(r,v)=\int_{\partial B_r} { |\nabla v |^2 \over |x|^{N-2}  },
\qquad J(r,v)=\int_{ B_r} { |\nabla v |^2 \over |x|^{N-2}  }.
$$
Then 
$$\varphi' ={\varphi \over r}
\left[-4 + {r I(r,(D_e u)^+) \over J(r,(D_e u)^+)}  + 
 {rI(r,(D_e u)^-) \over J(r,(D_e u)^-)}  \right]. \tag 7.2
$$
Now let $F$ be the fundamental solution of $(\Delta  + 2)  $. Then
$$F=F_1 |x|^{2-N} + F_2\log|x|,$$
 where $F_2$ is zero for odd $N$, and
$F_1$ and $F_2$ are regular functions not vanishing at the origin.
Next as in [ACF,  Lemma 5.1] (cf. also [Ca3]) we will have,
for a nonnegative subsolution $v$ to (7.1),
$$\align
J(r,v)& \leq 
(1+Cr)\int_{B(0,r)} ( \Delta(v^2/2) F +  v^2 F)\tag 7.3
\cr 
&=
(1+Cr) \int_{B(0,r)}( \Delta(v^2/2) F -  {v^2\over 2}\Delta F)\cr
&= (1+Cr) \int_{\partial B(0,r)} vv_n F - {v^2\over 2} F_n\cr
&\leq 
(1+Cr) \int_{\partial B(0,r)}( v|v_n| r^{2-N} + (N-2){v^2\over 2}
r^{1-N})\cr
&\leq (1+Cr){rI(r,v)\over2\gamma(\Gamma_v(r))},
\endalign 
$$
where $\gamma(\Gamma_v(r))$ is the corresponding $r$-dependent value
for $v$, as  in Lemma 2.2, and $v_n$ is the normal-directional
derivative of $v$ on the sphere $\partial B(0,r)$. 
Now let us consider the new function
$$\psi(r,D_e u) = \varphi (r, D_e u) e^{Cr},$$
with $C$ to be chosen later  so as to make $\psi$ nondecreasing.
Observe also that $(D_eu)^\pm$ are subsolutions to (7.1).
Differentiating we obtain
$$\psi'(r) = {\psi(r)\over r}
\left[ {r I(r,(D_e u)^+)\over J(r,(D_e u)^+)} +{r I(r,(D_e u)^-)\over
    J(r,(D_e u)^-)} -4 + Cr \right].\tag{7.4}
$$
Plugging (7.3) in (7.4) we obtain
$$\psi'\geq {\psi(r)\over r} \left[  2\gamma(\Gamma_+) (1 + O(r))+
  2\gamma(\Gamma_-) (1 + O(r)) -4 + Cr \right].$$
Now choosing $C$ large, we see that  the extra terms will be taken care of, and the
result follows as in [ACF, Lemma 5.1].

\AuthorRefNames [CyFrKiSi]
 
\vglue18pt

\references

[ACF]
\name{H.\ W. Alt, L.\ A. Caffarelli}, and \name{A. Friedman},
Variational problems with two phases and their free boundaries, 
{\it Trans.\ A.M.S.} {\bf 282} (1984), 431--461.

[AC]
\name{I. Athanasopoulos} and   \name{L.\ A. Caffarelli}, 
A theorem of real analysis and its application to free boundary problems,
{\it Comm.\ Pure Appl.\ Math.}\   {\bf  38 } (1985),  499--502.

[BKP]
\name{W. Beckner, C. Kenig}, and \name{J. Pipher},
in preparation.

[Ca1]
\name{L.\ A. Caffarelli},
The regularity of free boundaries in higher dimension,
{\it Acta Math}. {\bf 139} (1977), 155--184.

[Ca2]
\bibline,
Compactness methods in free boundary problems,
{\it Comm.\ P.D.E.} {\bf  5} (1980), 427--448.

[Ca3]
\bibline,
A Harnack inequality approach to the regularity of free boundaries,
{\it Scoula Norm.\ Sup.\ Pisa} {\bf 15  } (1988), 583--602.

[Ca4]
\bibline,
The obstacle problem revisited, III.  Existence theory,
compactness and dependence on $X$, 
{\it  J.\ Fourier Anal.\ Appl. }\  {\bf 4} (1998), 
383--402.

[CK]
\name{L.\ A. Caffarelli} and  \name{C. Kenig},
Gradient estimates for variable coefficient parabolic equations
and
singular perturbation problems,
{\it Amer.\ J.  of Math.}\  {\bf 120} (1998), 391--440.

[DF]
\name{E. di Benedetto}  and \name{A. Friedman},
Bubble growth in porous media,
{\it Indiana Univ.\ Math.\ J.} {\bf 35} (1986), 573--606.

[FH]
\name{S. Friedland} and  \name{W.\ K. Hayman},
Eigenvalue inequalities for the Dirichlet problem on spheres
and the growth of subharmonic functions,
{\it Comment.\ Math.\ Helv.}\  {\bf 51} (1976), 133--161.

[GT]
\name{D. Gilbarg} and \name{N.\ S. Trudinger},
{\it Elliptic Partial Differential Equations of Second Order}, $2^{\rm nd}$ edition,
 Springer-Verlag,  New York, 1983.

[GP]
 \name{B. Gutafsson} and  \name{M. Putinar},
The exponential transform and regularity of free boundaries in
two dimensions,
{\it Ann.\ Scoula Norm.\ Sup.\ Pisa Cl.\ Sci.}\ {\bf 26} (1998), 507--543.

[I1]
\name{V. Isakov},
{\it Inverse Source Problems{\rm ,} A.M.S.  Math.\ Surveys and Monographs} {\bf 34}, Providence, RI,
 1990.

[I2]
\bibline,
Inverse theorems on the smoothness of potentials,  in {\it Differential
Equations} {\bf 11} (1976), 50--57 (translated from Russian).

[KS1]
\name{L. Karp} and \name{H. Shahgholian},
Regularity of a free boundary problem,
 {\it J. Geometric Analysis},  to appear.

[KS2]
\bibline,
Regularity of a free boundary at the infinity point,
{\it Comm.\ Partial Differential Eq.}, to appear.

[KN]
\name{D. Kinderlehrer} and \name{L. Nirenberg},
Regularity in free boundary value problems, {\it Ann.\ Scuola\ Norm.\
Sup.\ Pisa Cl.\ Sci.}\ {\bf 4} (1977), 373--391.

[Ma]
\name{A.\ S. Margulis},
Potential theory for $L^p$-densities and its applications to inverse
problems of gravimetry, in
 {\it Theory and Practice of Gravitational and Magnetic Fields
Interpretation in USSR}, 188--197,
 Naukova Dumka Press, Kiev, 1983   (Russian).

[MY]
\name{J. McCarthy} and  \name{L. Yang},
Subnormal operators and quadrature domains, {\it Adv.\ Math.}\
{\bf 127 } (1997), 52--72.

[P]
 \name{M. Putinar},
Extremal solutions of the two dimensional $L$-problems
of moments, {\it J. Funct.\ Anal.}\
{\bf 136 } (1996), 331--364.

[Sa1]
\name{M. Sakai},
{\it Quadrature Domains{\rm ,}
Lecture Notes in Math}.\ {\bf 934}, Springer-Verlag, New York, 1982.

[Sa2]
 \bibline,
Regularity of a boundary having a Schwarz function, {\it Acta
Math.}\ {\bf 166 } (1991), 263--297.

[Sa3]
\bibline,
Regularity of boundaries of quadrature domains in two dimensions,
{\it SIAM J. Math.\ Anal.}\ {\bf 24 } (1993), 341--364.

[Shah]
 H. Shahgholian,
On quadrature domains and the Schwarz potential,
{\it J. Math.\ Anal.}\ {\bf 171} (1992), 61--78.

[Shap1]
 \name{H.\ S. Shapiro},
Global geometric aspects of Cauchy's  problem for the Laplace operator,
in {\it Geometrical  and Algebraical Aspects in Several Complex
Variables}, 309--324 , Cetraro,  1989. 

[Shap2]
\bibline,
{\it The Schwarz Function and Its Generalization to Higher Dimensions{\rm , }
Univ.\ of Arkansas Lecture Notes in the Math.\
Sci.}\ {\bf 9}, John Wiley \& Sons, Inc., New York, 1992.

[Spn]
\name{E. Sperner},
Zur symmetrisierung von funktionen auf sph\"aren,
{\it Math.\ Z.} {\bf 134} (1973), 317--327.

[Spk]
\name{J. Spruck},
Uniqueness in diffusion model of population biology,
{\it Comm.\ P.D.E.} {\bf 8} (1983), 1605--1620.

[St]
\name{V.\ N. Strakhov},
The inverse logarithmic potential problem for contact surface,
{\it Physics of the Solid Earth} {\bf 10} (1974), 104--114  (translated from
Russian).

[W]
\name{S. Williams},
Analyticity of the boundary for Lipschitz domains without
the Pompeiu property,
{\it Indiana Univ.\ Math.\ J}.
{\bf 30} (1981), 357--369.

\endreferences

\enddocument